\documentclass{rae}
\usepackage{amsmath}
\usepackage{amssymb}
\newcommand{\hk}{{\cal HK}}
\newcommand{\nbv}{{\cal NBV}}
\newcommand{\R}{{\mathbb R}}
\newcommand{\N}{{\mathbb N}}
\newcommand{\qed}{\mbox{$\quad\blacksquare$}}
\newtheorem{theorem}{Theorem}
\newtheorem{example}[theorem]{Example}
\newtheorem{remark}[theorem]{Remark}
\begin{document}
\begin{center}
{\large\bf A product convergence theorem for Henstock--Kurzweil integrals}
\vskip.25in
Parasar Mohanty\\
Erik Talvila\footnote{Research partially supported by the
Natural Sciences and Engineering Research Council of Canada.}\\ [2mm]
{\footnotesize
Department of Mathematical and Statistical Sciences\\
University of Alberta\\
Edmonton AB Canada T6G 2G1\\
pmohanty@math.ualberta.ca\\
etalvila@math.ualberta.ca}
\end{center}

{\footnotesize
\noindent
{\bf Abstract.}  Necessary and sufficient for $\int_a^bfg_n\to \int_a^bfg$ for
all Henstock--Kurzweil integrable functions $f$ is that $g$ be of
bounded variation, $g_n$ be uniformly bounded and of uniform bounded
variation and, on each compact interval in $(a,b)$, $g_n\to g$ in measure
or in the $L^1$ norm.
The same conditions are necessary and sufficient for $\|f(g_n-g)\|\to 0$
for all Henstock--Kurzweil integrable functions $f$.  If $g_n\to g$ a.e.
then convergence $\|fg_n\|\to\|fg\|$  for all Henstock--Kurzweil 
integrable functions $f$ is equivalent to
$\|f(g_n-g)\|\to 0$.
This extends a theorem due
to Lee Peng-Yee.
}\\
\vskip.25in

\noindent
Let $-\infty\leq a<b\leq\infty$ 
and denote the Henstock--Kurzweil integrable functions
on $(a,b)$ by $\hk$. 
The Alexiewicz norm of $f\in\hk$ is
$\|f\|=\sup_I|\int_I f|$ where the
supremum is taken over all intervals $I\subset(a,b)$.  If $g$ is a real-valued
function on $[a,b]$ we write $V_{[a,b]}g$ for the variation  of $g$ over 
$[a,b]$, dropping the subscript when the identity of $[a,b]$ is clear.  The
set of functions of normalised bounded variation, $\nbv$, 
consists of the functions on $[a,b]$ that are of bounded variation,
are left continuous and vanish at $a$.  
It is known that the multipliers  for $\hk$ are $\nbv$, i.e.,
$fg\in\hk$ for all $f\in\hk$ if and only if $g$ is equivalent to a
function in $\nbv$.
This paper is concerned with necessary and sufficient
conditions under which $\int_a^bfg_n\to\int_a^bfg$ for all $f\in\hk$.  One
such set of conditions was given by Lee Peng-Yee in 
\cite[Theorem~12.11]{lee}. If $g$ is of bounded variation, changing
$g$ on
a countable set will make it an element of $\nbv$.
With this
observation, a minor modification 
of Lee's theorem produces the following result.
\begin{theorem}\label{theorem1}{\rm \cite[Theorem~12.11]{lee}}
Let $-\infty<a<b<\infty$, let $g_n$ and $g$  be real-valued 
functions on $[a,b]$ with $g$ of bounded variation.
In order for $\int_a^bfg_n\to\int_a^bfg$ for all $f\in\hk$ it is necessary
and sufficient that
\begin{equation}
\left.
\begin{array}{l}
\text{for each interval } (c,d)\subset(a,b), \int_c^dg_n\to\int_c^dg 
\text{ as } n\to\infty,\\
\text{for each } n\geq 1, g_n \text{ is equivalent to a function } h_n\in\nbv,\\
\text{and there is  }M\in[0,\infty) \text{ such that } Vh_n\leq M 
\text{ for all } n\geq 1.
\end{array}
\right\}
\label{1}
\end{equation}
\end{theorem}

We extend this theorem to unbounded intervals, show that the condition
$\int_c^dg_n\to\int_c^dg$ in \eqref{1} can be replaced by $g_n\to g$ 
on each compact interval in $(a,b)$ either in 
measure or  in the $L^1$ norm, 
and that this also lets us conclude $\|f(g_n-g)\|\to 0$.
We also show that if $g_n\to g$ in measure or almost everywhere 
then $\|fg_n\|\to\|fg\|$  
for all $ f\in\hk$ if and only
if $\|fg_n-fg\|\to 0$ for all $f\in\hk$.

One might think the conditions \eqref{1} imply $g_n\to g$ almost
everywhere.  This is not the case, as is illustrated by the following
example \cite[p.~61]{folland}.
\begin{example}
{\rm
Let $g_n=\chi_{(j2^{-k},(j+1)2^{-k}]}$ where $0\leq j<2^k$ and $n=j+2^k$.
Note that $\|g_n\|_{\infty}=1$, $g_n\in\nbv$, $Vg_n\leq 2$, and 
$|\int_c^dg_n|\leq\|g_n\|=
2^{-k}< 2/n\to 0$, so that \eqref{1} is satisfied with $g=0$.  
For 
each $x\in(0,1]$ we
have $\inf_ng_n(x)=0$,  $\sup_ng_n(x)=1$, and  for no $x\in(0,1]$ does
$g_n(x)$ have a limit.  However, $g_n\to 0$ in measure since if
$T_n=\{x\in[0,1]: |g_n(x)|>\epsilon\}$ then for each  $0<\epsilon\leq 1$,
we have $\lambda(T_n)<2/n\to 0$ as $n\to \infty$ ($\lambda$ is Lebesgue
measure).
}
\end{example}

We have the following extension of Theorem~\ref{1}.
\begin{theorem}\label{theorem2}
Let $[a,b]$ be a compact interval in $\R$, let $g_n$ and $g$  be real-valued 
functions on $[a,b]$ with $g$ of bounded variation.
In order for $\int_a^bfg_n\to\int_a^bfg$ for all $f\in\hk$ it is necessary
and sufficient that
\begin{equation}
\left.
\begin{array}{l}
g_n\to g \text{ in measure as } 
n\to\infty,\\
\text{for each } n\geq 1\text{, }
g_n \text{ is equivalent to a function } h_n\in\nbv,\\
\text{and there is  }M\in[0,\infty) \text{ such that } Vh_n\leq M 
\text{ for all } n\geq 1.
\end{array}
\right\}
\label{2}
\end{equation}
If $(a,b)\subset\R$ is unbounded, then change the first line of
\eqref{2} by requiring 
$g_n\chi_I\to g\chi_I$  in measure for each compact interval $I\in (a,b)$. 
\end{theorem}

\bigskip
\noindent
{\bf Proof:} By working with $g_n-g$ we can assume $g=0$.  First consider
the case when $(a,b)$ is a bounded interval.

If $\int_a^b fg_n\to 0$ for all $f\in\hk$, then using Theorem~\ref{1} and
changing $g_n$ on a countable set, we can assume $g_n\in\nbv$, $Vg_n\leq M$,
$\|g_n\|_{\infty}\leq M$
and $\int_c^dg_n\to 0$ for each interval $(c,d)\subset(a,b)$.
Suppose $g_n$ does not converge to $0$ in measure.  Then there are $\delta,
\epsilon >0$ and an infinite index set ${\cal J}\subset\N$ such that
$\lambda(S_n)>\delta$ for each $n\in {\cal J}$, where $S_n=\{x\in(a,b):
g_n(x)>\epsilon\}$.  (Or else there is a corresponding set on which
$g_n(x)<-\epsilon$ for all $n\in{\cal J}$.)  Now let $n\in{\cal J}$.
Since $g_n$ is left continuous,
if $x\in S_n$ there is a number $c_{n,x}>0$ such that $[x-c_{n,x},x]\subset
S_n$.  Hence, $V_n:=\{[c,x]:x\in S_n \text{ and } [c,x]\subset S_n\}$ is a
Vitali cover of $S_n$.  So there is a finite set of 
disjoint closed intervals, $\sigma_n\subset V_n$, 
with $\lambda(S_n\setminus \cup_{I\in\sigma_n}I)
<\delta/2$.  Write $(a,b)\setminus \cup_{I\in\sigma_n}I=
\cup_{I\in\tau_n}I$ where $\tau_n$ is a set of disjoint open intervals with
${\rm card}(\tau_n)={\rm card}(\sigma_n)+1$.  Let $P_n={\rm card}(\{
I\in\tau_n:g_n(x)\leq\epsilon/2 \text{ for some } x\in I\})$.
Each interval $I\in\tau_n$ that does not have $a$ or $b$ as an endpoint
has contiguous intervals on its left and right that are in $\sigma_n$
(for each of which $g_n(x)>\epsilon/2$ for some $x$).
The interval $I$ then contributes more than $(\epsilon-\epsilon/2)+
(\epsilon-\epsilon/2)=\epsilon$ to the variation of $g_n$.  If $I$ has
$a$ as an endpoint then, since $g_n(a)=0$, $I$ contributes more than
$\epsilon$ to the variation of $g_n$.  If $I$ has
$b$ as an endpoint then $I$ contributes more than
$\epsilon/2$ to the variation of $g_n$.
Hence, 
$Vg_n\geq (P_n-1)
\epsilon+\epsilon/2=(P_n-1/2)\epsilon$.  (This inequality is still
valid if $P_n=1$.)  But, $Vg_n\leq M$ so $P_n\leq P$
for all $n\in{\cal J}$ and some $P\in\N$.  Then we have a set of intervals,
$U_n$, formed by taking unions of intervals from $\sigma_n$ and those
intervals in $\tau_n$ on which $g_n>\epsilon/2$.  Now, $\lambda(\cup_{I\in U_n
}I)>
\delta/2$, ${\rm card}(U_n)\leq P+1$ and $g_n>\epsilon/2$ on each
interval $I\in U_n$.  Therefore, there is an
interval $I_n\in U_n$ such that $\lambda(I_n)>\delta/[2(P+1)]$.  The
sequence of centres of intervals  $I_n$ has a convergent subsequence.  
There is then
an infinite index set ${\cal J}'\subset{\cal J}$ with the property that for
all $n\in{\cal J}'$ we have $g_n>\epsilon/2$ on an interval $I\subset(a,b)$
with $\lambda(I)>\delta/[3(P+1)]$.  Hence, $\limsup_{n\geq 1}\int_Ig_n>
\delta\epsilon/[6(P+1)]$.  This contradicts the fact that $\int_Ig_n\to 0$,
showing that indeed $g_n\to 0$ in measure.

Suppose \eqref{2} holds.  As above, we can assume $g_n\in\nbv$, $Vg_n\leq M$,
$\|g_n\|_{\infty}\leq M$ and $g_n\to 0$ in measure. Let $\epsilon>0$.
Define $T_n=\{x\in(a,b): |g_n(x)|>\epsilon\}$.  Then
\begin{eqnarray}
\left|\int_a^bg_n\right| & \leq & \int_{T_n}|g_n|+
\int_{(a,b)\setminus T_n}|g_n|\\
 & \leq & M\lambda(T_n)+\epsilon(b-a).
\end{eqnarray}
Since $\lim\lambda(T_n)=0$, it now follows that $\int_c^d g_n\to 0$ for
each $(c,d)\subset(a,b)$.  Theorem~\ref{theorem1} now shows $\int_a^b fg_n\to
0$ for all $f\in \hk$.

Now consider integrals on $\R$.  If $\int_{-\infty}^\infty fg_n\to 0$ for all
$f\in\hk$ then it is necessary that $\int_{a}^b fg_n\to 0$ for each
compact interval $[a,b]$.
By the current theorem,
$g_n\to g$  in measure on each $[a,b]$.  And, it is necessary that
$\int_{1}^\infty fg_n\to 0$.  The change of variables $x\mapsto 1/x$
now shows it is necessary that $g_n$ be equivalent to a function that
is uniformly bounded and of 
uniform bounded variation on $[1,\infty]$.  Similarly
with $\int_{-\infty}^1 fg_n\to 0$.  Hence, it is necessary that
$g_n$ be uniformly bounded and of uniform bounded variation on 
$\R$.

Suppose \eqref{2} holds with $g_n\to g$  in measure on each compact 
interval in $\R$. Write 
$\int_{-\infty}^\infty fg_n=\int_{-\infty}^a fg_n+\int_{a}^b fg_n+
\int_{b}^\infty fg_n$.  Use Lemma~24 in \cite{talvilafourier} to write
$|\int_{-\infty}^a fg_n|\leq \|f\chi_{(-\infty,a)}\|V_{[-\infty,a]}g_n\leq
\|f\chi_{(-\infty,a)}\|M\to 0$ as $a\to-\infty$.  We can then take a large
enough interval $[a,b]\subset\R$ and apply the current theorem on
$[a,b]$.
Other unbounded intervals are handled in a similar manner.
\qed

\begin{remark}\label{measure}
{\rm
If \eqref{2} holds then dominated convergence
shows $\|g_n-g\|_1\to 0$.  And, convergence in $\|\cdot\|_1$ implies 
convergence in measure.  Therefore, in the first statement of \eqref{2}
and in the last statement of Theorem~\ref{theorem2}, `convergence
in measure' can be replaced with `convergence in $\|\cdot\|_1$'.
Similar remarks apply to Theorem~\ref{theorem3}.
}
\end{remark}
\begin{remark}
{\rm
The change of variables argument in the second last paragraph of
Theorem~\ref{theorem2} can be replaced with an appeal to the Banach--Steinhaus
Theorem on unbounded intervals.  See \cite[Lemma~7]{sargent}.
Similarly in the proof of Theorem~\ref{theorem5}.
}
\end{remark}

The sequence of Heaviside step functions $g_n=\chi_{(n,\infty]}$ shows \eqref{2}
is not necessary to have $\int_{-\infty}^\infty fg_n\to 0$ for all
$f\in\hk$.  For then, $\int_{-\infty}^\infty fg_n=\int_n^\infty f\to 0$.
In this case, $g_n\in\nbv$ and $Vg_n=1$. 
However, $\lambda(T_n)=\infty$ for all $0<\epsilon<1$.  Note that
for each compact interval $[a,b]$ we have $\int_a^bg_n\to 0$ and $g_n\to 0$
in measure on $[a,b]$.

It is somewhat surprising that the conditions \eqref{2} are also necessary
and sufficient to have $\|f(g_n-g)\|\to 0$ for all $f\in\hk$.
\begin{theorem}\label{theorem3}
Let $[a,b]$ be a compact interval in $\R$, let $g_n$ and $g$  be real-valued
functions on $[a,b]$ with $g$ of bounded variation.
In order for $\|f(g_n-g)\|\to 0$ for all $f\in\hk$ it is necessary
and sufficient that
\begin{equation}
\left.
\begin{array}{l}
g_n\to g \text{ in measure as }
n\to\infty,\\
\text{for each } n\geq 1\text{, } 
g_n \text{ is equivalent to a function } h_n\in\nbv,\\
\text{and there is  }M\in[0,\infty) \text{ such that } Vh_n\leq M
\text{ for all } n\geq 1.
\end{array}
\right\}
\label{3}
\end{equation}
If $(a,b)\subset\R$ is unbounded, then change the first line of
\eqref{3} by requiring
$g_n\chi_I\to g\chi_I$  in measure for each compact interval $I\in (a,b)$.
\end{theorem}

\bigskip
\noindent
{\bf Proof:} Certainly \eqref{3} is necessary in order for
$\|f(g_n-g)\|\to 0$ for all $f\in\hk$.

If we have \eqref{3}, let $I_n$ be any sequence of intervals in
$(a,b)$.  We can again assume $g=0$.  Write ${\tilde g}_n=g_n\chi_{I_n}$.  Then
$\|{\tilde g}_n\|_{\infty}\leq\|g_n\|_{\infty}$, $V{\tilde g}_n\leq
Vg_n+2\|g_n\|_{\infty}$ and ${\tilde g}_n\to 0$ in measure.  
The result now follows by applying Theorem~\ref{theorem2} to $f{\tilde g}_n$.

Unbounded intervals are handled as in Theorem~\ref{theorem2}.
\qed

By combining Theorem \ref{theorem2} and Theorem \ref{theorem3} we 
have the following.
\begin{theorem}\label{theorem4}
Let $(a,b) \subset \R$ then $\int_a^b fg_n\to \int_a^b fg$ for all $f\in\hk$
if and only if $\|fg_n-fg\|\to 0$ for all $f\in\hk$.
\end{theorem}

Note that ${\displaystyle 
\|f(g_n-g)\|\geq \left|\,\|fg_n\|-\|fg\|\,\right|}$ so if $\|f(g_n-g)\|\to 0$
then $\|fg_n\|\to\|fg\|$.  Thus, \eqref{3} is sufficient to have
$\|fg_n\|\to
\|fg\|$ for all $f\in\hk$.  However, this condition is not necessary.
For example, let $[a,b]=[0,1]$. Define $g_n(x)=(-1)^n$.  Then 
$\|g_n\|_{\infty}=1$ and
$Vg_n=0$.  Let $g=g_1$.  
For no $x\in[-1,1]$ does the sequence $g_n(x)$ converge to $g(x)$.
For no open interval $I\subset[0,1]$ do we have
$\int_I(g_n-g)\to
0$.
And, $g_n$ does not converge to $g$ in measure.  However, let $f\in\hk$
with $\|f\|>0$.
Then $\|f(g_n-g)\|=0$
when $n$ is odd and when $n$ is even, $\|f(g_n-g)\|=2\|f\|$.
And yet, for all $n$, $\|fg_n\|=\|f\|=\|fg\|$.

It is natural to ask what extra condition should be given so that 
$\|fg_n\|\to\|fg\|$ will imply $\|fg_n-fg\|\to 0$.
We have the following.
\begin{theorem}\label{theorem5}
Let $g_n\to g$ in measure or almost everywhere.
Then $\|fg_n\|\to\|fg\|$  for all $ f\in\hk$ if and only 
if $\|fg_n-fg\|\to 0$ for all $f\in\hk$.
\end{theorem}
\noindent
{\bf Proof:} Let $[a,b]$ be a compact interval.
If $\|fg_n\|\to\|fg\|$  then $g$ is equivalent to $h\in\nbv$ 
\cite[Theorem~12.9]{lee} and for each $f\in\hk$ there is a constant
$C_f$ such that $\|fg_n\|\leq C_f$.  By the Banach--Steinhaus Theorem
\cite[Theorem~12.10]{lee},
each $g_n$ is equivalent to a function $h_n\in\nbv$ with $Vh_n\leq M$
and $\|h_n\|_\infty\leq M$.
Let $(c,d)\subset(a,b)$.  
By dominated convergence, $\int_c^dg_n\to \int_c^dg$.  It now follows from
Theorem~\ref{theorem1} that $\int_a^bfg_n\to\int_a^bfg$ for all
$f\in\hk$.  Hence, by Theorem~\ref{theorem4}, $\|fg_n-fg\|\to 0$
for all $f\in\hk$.

Now suppose $(a,b)=\R$ and $\|fg_n\|\to\|fg\|$ for all $f\in\hk$.
The change of variables $x\mapsto 1/x$ shows the Banach--Steinhaus
Theorem still holds on $\R$.  We then have each $g_n$ equivalent
to $h_n\in\nbv$ with $Vh_n\leq M$ and $\|h_n\|_\infty\leq M$.  
As with the end of the proof of
Theorem~\ref{theorem2}, given $\epsilon >0$ we can find $c\in\R$
such that $|\int_{-\infty}^c fg_n|<\epsilon$ for all $n\geq 1$.
The other cases are similar.
\qed\\

\noindent
{\bf Acknowledgment.} An anonymous referee provided reference \cite{sargent}
and pointed out that in place
of convergence in measure we can use convergence in $\|\cdot\|_1$
(cf. Remark~\ref{measure}).

\end{document}